\documentclass[12pt]{article}
\usepackage{amsfonts,amssymb}

\newtheorem{theorem}{Theorem}[section]
\newtheorem{prop}[theorem]{Proposition}

\newtheorem{lemma}[theorem]{Lemma}

\newcommand{\qed}{\hfill$\Box$}

\newtheorem{unnumber}{}

\newtheorem{prepf}[unnumber]{Proof}
\newenvironment{pf}{\prepf\rm}{\qed\endprepf}

\newcommand{\Aut}{\mathop{\mathrm{Aut}}}
\newcommand{\sd}{\mathbin{\triangle}}
\newcommand{\Sym}{\mathop{\mathrm{Sym}}}
\newcommand{\Alt}{\mathop{\mathrm{Alt}}}

\newcommand{\Fix}{\mathop{\mathrm{Fix}}\nolimits}

\newcommand{\AGL}{\mathop{\textrm{AGL}}}

\renewcommand{\wr}{\mathop{\textrm{wr}}}

\newcommand{\fix}{\mathop{\mathrm{fix}}\nolimits}
\newcommand{\orb}{\mathop{\mathrm{orb}}\nolimits}

\def\cent#1#2{{\bf C}_{{#1}}{{(#2)}}}

\begin{document}

\title{Most switching classes with primitive automorphism groups\\
contain graphs with trivial groups}
\author{Peter J. Cameron and Pablo Spiga\\[5mm]
School of Mathematics and Statistics\\University of
St Andrews\\St Andrews, Fife  KY16 9SS\\Scotland\\[3mm]
Dipartimento di Matematica e Applicazioni\\University of
Milano-Bicocca\\Milano, 20125 Via Cozzi 55\\Italy\\[3mm]}
\date{To the memory of \'Akos Seress}
\maketitle

\begin{abstract}
The operation of \emph{switching} a graph $\Gamma$ with respect to a subset
$X$ of the vertex set interchanges edges and non-edges between $X$ and its
complement, leaving the rest of the graph unchanged. This is an equivalence
relation on the set of graphs on a given vertex set, so we can talk about the
automorphism group of a \emph{switching class} of graphs.

It might be thought that switching classes with many automorphisms would have
the property that all their graphs also have many automorphisms. However the
main theorem of this paper shows a different picture: with finitely many
exceptions, if a non-trivial switching class $\mathcal{S}$ has primitive
automorphism group, then it contains a graph whose automorphism group is
trivial. We also find all the exceptional switching classes; up to 
complementation, there are just six.
\end{abstract}

\section{Introduction}

The purpose of this paper is to prove the following theorem:

\begin{theorem}\label{thrm:main}
Let $\mathcal{S}$ be a switching class of graphs on $n$ vertices, whose
automorphism group is primitive but not the symmetric group. Then, with
just six exceptions (up to complementation), there is a graph
$\Gamma\in\mathcal{S}$ with trivial automorphism group.

The six exceptions are as follows (we give the degree and automorphism group
of the switching class in each case): $(5,D_{10})$, $(6,\mathrm{PSL}(2,5))$,
$(9,S_3\wr S_2)$, $(10, \mathrm{P}\Sigma\mathrm{L}(2,9)\cong S_6)$,
$(14,\mathrm{PSL}(2,13))$, and $(16, 2^4.S_6\cong2^4.\mathrm{Sp}(4,2))$.
\end{theorem}

The next two sections explain the context and give definitions. Then we
prove that the list of exceptions is finite, and lastly we find all the
exceptions.

Note that the assumption of primitivity is necessary: the switching class
containing the complete multipartite graph $K_{r,r,\ldots,r}$ with $r\ge3$
contains no graph with trivial automorphism group. (The switching set $X$
or its complement meets a set of the multipartition in at least two points,
and the permutation transposing two such points is an automorphism.)

\section{Primitive permutation groups}

One of the most useful consequences of the Classification of the Finite Simple
Groups (CFSG) is that primitive permutation groups (other than the symmetric
and alternating groups) are ``small''. In fact, elementary combinatorial
arguments due to Babai~\cite{babai} for uniprimitive groups and
Pyber~\cite{pyber} for multiply transitive groups, provide bounds which are
good enough for many purposes; but the results can be strengthened by
using CFSG. The best result, which we use later, is due to
Mar\'oti~\cite{maroti}. (Babai's bound asserts that if $G$ is primitive but
not $2$-transitive of degree $n$, then $|G|\le n^{4\sqrt{n}\log_2 (n)}$. This is
best possible apart from the logarithm in the exponent, but Mar\'oti's result
gives a bound of $n^{1+\log_2 (n)}$ with known exceptions.) We state Mar\'oti's
theorem later.

A \emph{permutation group} is a subgroup of the symmetric group $S_n$ of
degree $n$, acting on a domain $\Omega$ of size $n$. It is \emph{transitive}
if it leaves no subset of $\Omega$ invariant except for $\emptyset$ and
$\Omega$, and \emph{primitive} if, in addition, it leaves no partition of
$\Omega$ invariant apart from the partition into singletons and the partition
into a single part.

The line of research reported here began with a theorem of Cameron, Neumann
and Saxl~\cite{cns}:

\begin{theorem}
If $G$ is a primitive group on $\Omega$, other than $S_n$ and
$A_n$ and finitely many exceptions, then there is a subset of $\Omega$ whose
setwise stabiliser in $G$ is the identity.
\end{theorem}

This result has been quantified in various ways:
\begin{enumerate}\itemsep0pt
\item Seress~\cite{seress} found all the exceptions: there are $43$ of them,
the largest degree being $32$. (We give a slight strengthening of this result
in Section~4.)
\item Cameron~\cite{regorbs} showed that the proportion of subsets whose
stabiliser is the identity tends to $1$ as $n\to\infty$ (in primitive groups of
degree $n$ other than $S_n$ and $A_n$).
\item Babai and Cameron~\cite{bc} showed that we can take the size of the
subset to be at most $n^{1/2+o(1)}$.
\end{enumerate}

In his paper, Seress gives the numbers of the primitive groups in the
\textsf{GAP} computer algebra system. However, the numbers have changed. For
the reader's convenience, we repeat the list, with the numbers in the current
version 4.7.4 of \textsf{GAP}. Each entry has the form $(n,m,G)$, where
$G$ is \texttt{PrimitiveGroup}$(n,m)$ in \textsf{GAP}.
\begin{quote}
$(5,2,D_{10})$,
$(5,3,\mathrm{AGL}(1,5))$,
$(6,1,\mathrm{PSL}(2,5))$,
$(6,2,\mathrm{PGL}(2,5))$,
$(7,4,\mathrm{AGL}(1,7))$,
$(7,5,\mathrm{PSL}(3,2))$,
$(8,2,\mathrm{A}\Gamma\mathrm{L}(1,8))$,
$(8,4,\mathrm{PSL}(2,7))$,
$(8,5,\mathrm{PGL}(2,7))$,
$(8,3,\mathrm{AGL}(3,2))$,
$(9,2,3^2.D_8=S_3\wr S_2)$,
$(9,5,\mathrm{A}\Gamma\mathrm{L}(1,9))$,
$(9,6,\mathrm{ASL}(2,3))$,
$(9,7,\mathrm{AGL}(2,3))$,
$(9,8,\mathrm{PSL}(2,8))$,
$(9,9,\mathrm{P}\Gamma\mathrm{L}(2,8))$,
$(10,2,S_5)$,
$(10,3,\mathrm{PSL}(2,9))$,
$(10,5,\mathrm{P}\Sigma\mathrm{L}(2,9))$,
$(10,4,\mathrm{PGL}(2,9))$,
$(10,6,M_{10})$,
$(10,7,\mathrm{P}\Gamma\mathrm{L}(2,9))$,
$(11,1,\mathrm{PSL}(2,11))$,
$(11,2,M_{11})$,
$(12,4,\mathrm{PGL}(2,11))$,
$(12,1,M_{11})$,
$(12,2,M_{12})$,
$(13,7,\mathrm{PSL}(3,3))$,
$(14,2,\mathrm{PGL}(2,13))$,
$(15,4,\mathrm{PSL}(4,2))$,
$(16,12,\mathrm{A}\Gamma\mathrm{L}(2,4))$,
$(16,17,2^4.A_6)$,
$(16,16,2^4.S_6)$,
$(16,20,2^4.A_7)$,
$(16,11,\mathrm{AGL}(4,2))$,
$(17,7,\mathrm{PSL}(2,16).2)$,
$(17,8,\mathrm{P}\Gamma\mathrm{L}(2,16))$,
$(21,7,\mathrm{P}\Gamma\mathrm{L}(3,4))$,
$(22,1,M_{22})$,
$(22,2,M_{22}.2)$,
$(23,5,M_{23})$,
$(24,1,M_{24})$,
$(32,3,\mathrm{AGL}(5,2))$.
\end{quote}

\section{Switching}

This section gives a very brief introduction to graph switching and some of its
many applications. See Seidel~\cite{seidel} or Taylor~\cite{taylor} for more
details, including connections with equiangular lines, group cohomology, graph
spectra, and finite simple groups.

Switching of graphs was introduced by Seidel in connection with equiangular
lines in Euclidean space, and later used in his classification of strongly
regular graphs with least eigenvalue $-2$. At about the same time, Higman
introduced an equivalent concept, that of a ``two-graph'' (which we describe
below), to give a direct construction of the Conway group $Co_3$

Let $\Gamma$ be a graph on the vertex set $\Omega$. For a subset $X$ of
$\Omega$, the operation $\sigma_X$ of \emph{switching} with respect to $X$
replaces edges between $X$ and $\Omega\setminus X$ by non-edges and non-edges
by edges, leaving edges and non-edges within or outside $X$ unchanged. Since
$\sigma_X\sigma_Y=\sigma_{X\sd Y}$, switching is an equivalence relation on
the set of all graphs on the vertex set $\Omega$. Since
$\sigma_X=\sigma_{\Omega\setminus X}$, there are $2^{n-1}$ graphs in a switching
class, where $n=|\Omega|$.

\paragraph{Example}
A regular icosahedron has six diagonals, any two making the same
angle. Choose one of the two vertices on each diagonal. The corresponding
induced subgraph of the skeleton of the icosahedron has one of four possible
forms: a pentagon with an isolated vertex; a pentagon with a vertex joined
to all others; a triangle with pendant edges at each vertex; and the 
complement of the preceding. These graphs form a switching class; the numbers
of graphs of each type are $6$, $6$, $10$, $10$ respectively. (For switching
with respect to a subset $X$ of the six diagonals corresponds to changing 
the choice of vertex on the diagonals in $X$.) Any switching class can be
represented by a set of equiangular lines in this way; see Seidel~\cite{seidel}
for a description of this.

\medskip

Given a graph $\Gamma$ with vertex set $\Omega$, the $3$-uniform hypergraph on
$\Omega$ whose hyperedges are the $3$-subsets of $\Omega$ containing  an odd
number of edges of $\Gamma$ is unaffected by switching; indeed, two
graphs belong to the same switching class if and only if they give the same
$3$-uniform hypergraph in this way. The hypergraphs which arise are
characterised by the property that any $4$-subset of $\Omega$ contains an even
number of hyperedges. Graham Higman, who introduced this notion, called such an
object a \emph{two-graph}. He gave a combinatorial construction and uniqueness
proof for a two-graph on $276$ points whose automorphism group is Conway's
third group $Co_3$. (The corresponding set of $276$ equiangular lines occurs in
the Leech lattice.)

The \emph{automorphism group} of a switching class consists of all permutations
of the vertex set which map the graphs of the class among themselves. It
suffices to assume that a single graph in the switching class is mapped to
a graph in the class. The automorphism group of a switching class coincides
with the automorphism group of the associated two-graph (as $3$-uniform 
hypergraph. The automorphism group of a switching class contains the
automorphism groups of all the graphs in the class as subgroups. Indeed, there
is a Frucht-style \emph{representation theorem} for a group and its
subgroups~\cite{c:aut}:

\begin{prop}
Given any
finite group $G$, there is a switching class $\mathcal{S}$ with the
properties:
\begin{enumerate}\itemsep0pt
\item $\Aut(\mathcal{S})\cong G$;
\item identifying these two groups, for any subgroup $H\le G$, there is a graph
$\Gamma\in\mathcal{S}$ with $\Aut(\Gamma)=H$.
\end{enumerate}
\end{prop}

Switching classes can be divided into two types, which we will call Type~I and
Type~II; the types are distinguished by the vanishing or non-vanishing of
a certain cohomology class (see Cameron~\cite{coho}).

\paragraph{Type I:} A Type~I switching class $\mathcal{S}$ contains a graph
$\Gamma$ such that $\Aut(\Gamma)=\Aut(\mathcal{S})$. For any subset $X$ of the
vertex set $V$ of $\Gamma$, $\Aut(\sigma_X(\Gamma))$ is the stabiliser of the
partition $\{X,V\setminus X\}$ in $\Aut(\Gamma)$. Every switching class on
an odd number of vertices is of Type~I: for such a switching class contains
a unique graph $\Gamma$ in which all vertices have even valency, and the
partition $\{X,V\setminus X\}$ is the partition of $V$ into vertices of
even and odd valency in $\sigma_X(\Gamma)$. The switching classes of
the $5$-cycle and the $3\times3$ grid (the line graph of $K_{3,3}$) are
of this type, and provide two of the examples in Theorem~\ref{thrm:main}; their
automorphism groups are $D_{10}$ and $S_3\wr S_2$ respectively.

\paragraph{Type II:} The automorphism group of a Type~II switching class
$\mathcal{S}$ properly contains the automorphism group of any graph in
$\mathcal{S}$. It can happen that the automorphism group of a switching
class is $2$-transitive; such a class, if not trivial (i.e. the switching
class of the complete or null graph), is necessarily of Type~II. Note,
however, that no non-trivial switching class can have a $3$-transitive
automorphism group, since $\Aut(\mathcal{S})$ is equal to the automorphism
group of the corresponding two-graph.

\smallskip

We now briefly describe the other four examples occurring in our theorem.
(The examples on $5$ and $9$ points are of Type~I and are described above.)
\begin{itemize}
\item Let $q$ be a prime power congruent to $1$~(mod~$4$). The group
$\mathrm{P}\Sigma\mathrm{L}(2,q)$, acting on $q+1$ points, has two orbits
on $3$-subsets; each orbit forms a two-graph, which is self-complementary
(that is, the two-graphs are isomorphic). The examples with $6$, $10$ and
$14$ points are of this type. The example on $6$ points corresponds to the
six diagonals of the icosahedron, as described earlier. The example on $10$
points corresponds to the switching class containing the Petersen graph (or
its complement, the line graph of $K_5$).
\item The example on $16$ points is a \emph{symplectic} two-graph. The
points are the vectors of a $4$-dimensional vector space over $\mathbb{F}_2$,
carrying a symplectic form (a non-degenerate alternating bilinear form) $B$.
The triples of the two-graph are those $\{x,y,z\}$ for which
\[B(x,y)+B(y,z)+B(z,x)=0.\]
It can also be described as the switching class of the Clebsch graph. Its
automorphism group is the semidirect product of the translation group of
$V$ by the symplectic group $\mathrm{Sp}(4,2)\cong S_6$.
\end{itemize}

\section{Bounding the exceptions}

We begin the proof of Theorem~\ref{thrm:main}. In this section we show that 
exceptions have degree at most $32$.

Let $\Omega$ be a finite set and let $g$ be a permutation on $\Omega$. We 
denote by $\Fix_\Omega(g)$ the set $\{\omega\in \Omega\mid \omega^g=\omega\}$,
by $\fix_\Omega(g)$ the cardinality $|\Fix_\Omega(g)|$ and  by $\orb_\Omega(g)$
the number of cycles of $g$   (in its decomposition in disjoint cycles). 

Given a subset $X$ of $\Omega$, we denote by $G_X$ the set-wise stabiliser
$\{g\in G\mid X^g=X\}$. We denote by $2^\Omega$ the power-set of $\Omega$,
that is, the set of subsets of $\Omega$. Moreover, we define 
\[
\mathcal{F}(G)=\{X\in 2^\Omega\mid X^g=X \textrm{ for some }g\in G\setminus\{1\}\}.
\]

\begin{lemma}\label{basic3}
Let $g$ be a permutation of $\Omega$ and let $p$ the smallest prime dividing
the order of $g$. Then
\[\orb_\Omega (g)\leq \frac{|\Omega|}{p}+\frac{p-1}{p}\fix_\Omega(g)\]
and
\[2^{\orb_\Omega(g)}=|\{X\in 2^\Omega\mid X^g=X\}|.\]
In particular,
\[|\{X\in 2^\Omega\mid X^g=X\}|\leq 2^{\frac{|\Omega|}{p}+\frac{p-1}{p}\fix_\Omega(g)}\]
and 
\[
|\mathcal{F}(G)|\leq \sum_{g\in G\setminus\{1\}}2^{\orb_\Omega(g)}\leq \sum_{g\in G\setminus\{1\}}2^{\frac{|\Omega|}{p}+\frac{p-1}{p}\mathrm{fix}_\Omega(g)}.
\]
\end{lemma}

\begin{pf}
The element $g$ has cycles of size $1$ on $\Fix_\Omega(g)$ and of size at
least $p$ on $\Omega\setminus \Fix_\Omega(g)$. Thus
\[\orb_\Omega(g)\leq |\Fix_\Omega(g)|+\frac{|\Omega\setminus\Fix_\Omega(g)|}{p}
=\frac{|\Omega|}{p}+\frac{p-1}{p}\fix_\Omega(g).\]
The rest of the lemma is obvious.
\end{pf}

We require a slight strengthening of the result of Seress~\cite{seress}
on primitive groups with no regular orbit on the power set.

\begin{prop}
\label{Akos}
Let $G$ be a finite primitive group on $\Omega$ with $\Alt(\Omega)\nleq G$.
Then either there exists a subset $X$ of $\Omega$ with $|X|< n/2$ 
and $G_X=1$, or $G$ is one of the groups listed
in~{\rm\cite[Theorem~$2$]{seress}}, or
$|\Omega|=16$ and $G=2^4.\mathrm{SO}_4^{-}(2)$.
\end{prop}

\begin{pf}
Let $n$ be the cardinality of $\Omega$. Observe that if there exists
 $X\in 2^\Omega$ with $G_X=1$, then replacing $X$ by $\Omega\setminus X$ if
necessary we have $|X|\leq n/2$ and we still have $G_X=1$.

If $n$ is odd, then the proof follows from~\cite[Theorem~$2$]{seress}.
Suppose then that $n$ is even.  Let $M$ be a subgroup of $\Sym(\Omega)$
maximal (with respect to set-inclusion) such that $G\leq M$ and
$M\notin\{\Alt(\Omega),\Sym(\Omega)\}$. Observe that if $M_X=1$ for some $X\in 2^\Omega$, then also $G_X=1$. Replacing $G$ by $M$ if necessary,
we assume that $G=M$. Now, the structure of $M$ is well understood
(see for example~\cite{LPS1}):
\begin{enumerate}
\item $n=m^r$ for some positive integers $m$ and $r$ with $m\geq 5$ and
$r\geq 2$, and $G$  is permutation isomorphic to $S_m\wr S_r$ with
its natural product action;
\item $G$ is the stabiliser of a diagonal structure, that is, $G$ is a group
of ``diagonal type'';
\item $G$ is the stabiliser of an affine structure, that is, $n=2^d$ for some positive integer $d$, and $G$ is permutation isomorphic to
$\AGL_d(2)$ endowed of its natural ``affine'' action;
\item $G$ is almost simple.
\end{enumerate}
We deal with each of these cases in turn.

Suppose that Case~(a) holds. Then from the proof of~\cite[Lemma~$4$]{seress}
we see that there exists $X\in2^\Omega$ with $G_X=1$ and with
$|X|=4m-5$ when $r=2$, and
$|X|=\sum_{i=0}^{r-1}(m-2)^im^{r-1-i}+3m-5$ when $r\geq 3$. An easy
calculation gives $|X|<n/2$ unless $(m,r)=(6,2)$. When $m=6$, we have $|X|=19$, $|\Omega\setminus X|=11<n/2$ and 
$G_{\Omega\setminus X}=1$.

Suppose that Case~(c) holds. Then $G\cong \AGL(d,2)$ for some positive integer
$d$. Now, a non-identity element of $G$ fixes at most $n/2$ points and
hence fixes at most $2^{3n/4}$ elements of $2^\Omega$ by
Lemma~\ref{basic3}. Therefore $|\mathcal{F}(G)|\leq 2^{3n/4}|G|$. Now a 
computation gives $2^{3n/4}|G|<2^n-{n\choose n/2}$ when $d\geq 9$. In
particular, for $d\geq 9$, there exists $X\in2^\Omega$ with
$|X|<n/2$ and $G_X=1$. When $d\leq 8$, the proof follows by
computation.

Suppose that Case~(d) holds. Assume that $G=S_m$ in its action on the
$r$-subsets of $\{1,\ldots,m\}$, with $2\leq r<m/2$. Then
Seress~\cite[Lemma~9]{seress} shows that there exists $X\in2^\Omega$
such that $G_X=1$ and with $|X|=m-r+1$ when $r\in \{2,3\}$ and
$|X|\leq 2(m-r+1)$ when $r\geq 4$. In both cases
$|X|<\frac{1}{2}{m\choose r}=n/2$. 

When $G=M_{n}$ (where $M_n$ is the Mathieu group of degree $n$)  the proof
follows from a computation.

Suppose that Case~(b) holds, or that Case~(d) holds and $G$ is neither
$S_m$ in its action on the $r$-subsets of $\{1,\ldots,m\}$ nor the Mathieu
group of degree $n$. Now from~\cite[Corollary~$1$]{GM} we have
$\fix_\Omega(g)\leq 4n/7$ for every $g\in G\setminus\{1\}$, and 
from~\cite{maroti} we have $|G|\leq n^{1+\log_2(n)}$.  Thus
Lemma~\ref{basic3} gives $|\mathcal{F}(G)|\leq 2^{11n/14}n^{1+\log_2(n)}$.
Now a computation gives that $2^{11n/14}n^{1+\log_2(n)}<2^n-{n\choose n/2}$
when $n\geq 386$. In particular, for $n\geq 386$, there exists
$X\in2^\Omega$ with $X\notin \mathcal{F}(G)$ and
$|X|\neq n/2$, and the proof follows immediately. 

Finally the primitive groups of diagonal type and almost simple of degree
less than $386$ can be easily checked with a computer.
\end{pf}

Observe that $\mathrm{SO}_4^{-}(2)\cong S_5$. There are two non-isomorphic
primitive groups of degree $16$ and with point stabiliser isomorphic to
$S_5$: one $2$-transitive and the other (namely
$2^4.\mathrm{SO}_4^{-}(2)$ in the statement of Proposition~\ref{Akos})
having rank $3$. The latter is \texttt{PrimitiveGroup}$(16,18)$ in the current
\textsf{GAP} list.

\medskip

We now begin the proof of the main theorem. We will show that an exception has
degree at most $32$.

In the following proof, with a slight abuse of terminology, we say that $M$
is a {\em maximal subgroup} of $\Sym(\Omega)$ if $M\leq\Sym(\Omega)$,
$\Alt(\Omega)\nleq M$ and either $\Alt(\Omega)$ or $\Sym(\Omega)$ is the
only subgroup of $\Sym(\Omega)$ containing $M$. 

\begin{pf}{\bf of Theorem~\ref{thrm:main}: the degree is at most $32$. }
Let $\mathcal{S}$ be a switching class of graphs on $n$ vertices, whose automorphism group $G$ is primitive but not the symmetric group. Let $\Omega$ be the domain of $G$, let $\mathbb{F}_2$ be the finite
field of cardinality $2$ and let $V$ be the permutation $\mathbb{F}_2G$-module for the
action of $G$ on $\Omega$. Thus $V$ has basis $(e_\omega\mid \omega\in \Omega)$
indexed by the elements of $\Omega$. Let $e=\sum_{\omega\in \Omega}e_\omega$
and let $W$  be the quotient $G$-module $V/\langle e\rangle$. 

\smallskip

The proof is divided into the two parts we called Type~I and Type~II earlier;
recall that a switching class has Type~I if its automorphism group is equal
to the automorphism group of some graph in the class.

\paragraph{Type I:} Thus $G=\Aut(\mathcal{S})=\Aut(\Gamma_0)$, for
some $\Gamma_0\in\mathcal{S}$.

If there exists $X\in 2^\Omega$ with
$|X|<n/2$ and $G_X=1$, then the proof follows immediately: the 
graph $\sigma_X(\Gamma_0)$ lies in $\mathcal{S}$ and has trivial
automorphism group. In particular, in view of Proposition~\ref{Akos} we may
assume that $G$ is one of the groups listed in~\cite[Theorem~$2$]{seress},
or $|\Omega|=16$ and $G= 2^4.\mathrm{SO}_4^-(2)$. 

A quick inspection reveals that the groups in~\cite[Theorem~$2$]{seress} have
degree at most $32$.

\paragraph{Type II:}
Assume that no graph in $\mathcal{S}$ has trivial automorphism group. Thus
$G$ has no regular orbit on $\mathcal{S}$ and hence
\begin{equation}\label{hopefullylast}
\#G\textrm{-orbits on }\mathcal{S}\geq 2\cdot \frac{|\mathcal{S}|}{|G|}=\frac{2^n}{|G|}.
\end{equation}

A result of Mallows and Sloane~\cite{ms} asserts that any permutation $g$
which fixes a switching class $\mathcal{S}$ fixes a graph
$\Gamma\in \mathcal{S}$. Moreover, another
graph, $\Gamma'\in \mathcal{S}$ is fixed by $g$ if and only if the partition
$\{X,\Omega\setminus X\}$ which switches $\Gamma$ to $\Gamma'$ is
fixed by $g$. It follows that the number of graphs in $\mathcal{S}$ fixed by
$g$ is equal to the number of vectors of $W$ fixed by $g$; that is,
\[\mathrm{fix}_{\mathcal{S}}(g)=\mathrm{fix}_W(g)\quad\textrm{for every }g\in G.\]
It follows that the permutation characters of $G=\Aut(\mathcal{S})$ on the
switching class $\mathcal{S}$ and on the ``switching module'' $W$ are equal.
Therefore 
\begin{equation}\label{eq:eq0}
\#G\textrm{-orbits on }\mathcal{S}=\# G\textrm{-orbits on }W.
\end{equation}

The definition of $W$ gives 
\begin{eqnarray}\label{eq:eq1}\nonumber
\#G\textrm{-orbits on }W&=&\frac{1}{2}\left(\right.\#G\textrm{-orbits on }2^\Omega\\
&&\qquad\left.+\#\textrm{self-complementary }G\textrm{-orbits on }2^\Omega\right),
\end{eqnarray}
where a $G$-orbit is self-complementary if it contains the complement of each
of its elements. From the Orbit-Stabiliser lemma and Lemma~\ref{basic3}, the
first summand in Equation~(\ref{eq:eq1}) is
$\sum_{g\in G}2^{\orb_\Omega(g)}/|G|$. 

A subset $X$ of $\Omega$  lies in a self-complementary orbit if and only
if there is a fixed-point-free element $g\in G$ of $2$-power order with
$\Omega=X\cup X^g$ and $X\cap X^g=\emptyset$. The
number of such sets $X$ arising from a given element $g$ is
$2^{\orb_\Omega(g)}$. Denote by $G_2$ the set
\[\{g\in G\mid |g|=2^\ell \textrm{ for some }\ell>1 \textrm{ and }\fix_\Omega(g)=0\}.\]
Thus, the second summand in Equation~(\ref{eq:eq1}) is at most
$\sum_{g\in G_2}2^{\orb_\Omega(g)}$.

Now from Equations~(\ref{hopefullylast}),~(\ref{eq:eq0}) and~(\ref{eq:eq1})
we obtain
\begin{equation}\label{eq:eq2}
\frac{2^n}{|G|}\leq \frac{1}{2}\left(\frac{1}{|G|}\sum_{g\in G}2^{\orb_\Omega(g)}+\sum_{g\in G_2}2^{\orb_\Omega(g)}\right).
\end{equation} 

Before dealing with the general case we deal separately with a few special
situations. The case division here is dictated by the theorem of
Mar\'oti~\cite{maroti} on the orders of primitive groups. For convenience we
state Mar\'oti's result here.

\begin{theorem}
Let $G$ be a primitive permutation group of degree~$n$. Then one of the
following holds:
\begin{enumerate}
\item $G$ is a subgroup of $S_m\wr S_r$ containing $(A_m)^r$, where the action
of $S_m$ is on $k$-subsets of $\{1,\ldots,m\}$ and the wreath product
has the product action of degree $n={m\choose k}^r$;
\item $G=M_{11}$, $M_{12}$, $M_{23}$ or $M_{24}$ with their $4$-transitive
action;
\item $\displaystyle{|G|\le n\cdot\prod_{i=0}^{\lfloor\log_2(n)\rfloor-1}
(n-2^i)<n^{1+\lfloor\log_2(n)\rfloor}}$.
\end{enumerate}
\end{theorem}

\subparagraph{Case (a)(i):} the socle of $G$ is $A_m$ in its action on
the $k$-subsets of $\{1,\ldots,m\}$ with $2\leq k< m/2$ and $m\geq 5$. 

\smallskip

Suppose first that $k=2$. Hence $n={m\choose 2}$ and we may identify
$\Omega$ with the set of $2$-subsets of $\{1,\ldots,m\}$. Let $\tau$ be the 
two-graph induced by the switching class $\mathcal{S}$. 
Suppose that $m\geq 9$. Now, $G$ has rank $3$ in its action on $\Omega$ and
$G$ has five orbits on the set of $3$-subsets of $\Omega$: namely
\begin{eqnarray*}
\mathcal{O}_1&=&\{\{a,b\},\{a,c\},\{a,d\}\mid a,b,c,d\textrm{ distinct elements of }\{1,\ldots,m\}\},\\
\mathcal{O}_2&=&\{\{a,b\},\{b,c\},\{a,c\}\mid a,b,c\textrm{ distinct elements of }\{1,\ldots,m\}\},\\
\mathcal{O}_3&=&\{\{a,b\},\{b,c\},\{c,d\}\mid a,b,c,d\textrm{ distinct elements of }\{1,\ldots,m\}\},\\
\mathcal{O}_4&=&\{\{a,b\},\{a,c\},\{d,e\}\mid a,b,c,d,e\textrm{ distinct elements of }\{1,\ldots,m\}\},\\
\mathcal{O}_5&=&\{\{a,b\},\{c,d\},\{e,f\}\mid a,b,c,d,e,f\textrm{ distinct elements of }\{1,\ldots,m\}\}.
\end{eqnarray*} 
(These orbits correspond to the  non-isomorphic graphs with three edges.)
Replacing $\tau$ by its complement 
$\tau'=\{x\subseteq\Omega\mid |x|=3, x\notin \tau\}$, 
we may assume  that $\mathcal{O}_1\subseteq \tau$. Now consider the $4$-subset
$x=\{\{a,b\},\{a,c\},\{a,d\},\{e,f\}\}$ of $\Omega$. Clearly, $x$ has one
$3$-subset in common with $\mathcal{O}_1$, zero with $\mathcal{O}_2$,
$\mathcal{O}_3$ and $\mathcal{O}_5$, and three with $\mathcal{O}_4$. As
$\tau$ is a two-graph, $x$ has an even number of $3$-subsets in common with
$\tau$ and hence $\mathcal{O}_4\subseteq \tau$. An entirely similar argument
using $x=\{\{a,b\},\{a,c\},\{a,d\},\{b,c\}\}$ yields
$\mathcal{O}_2\subseteq \tau$. Now an easy inspection gives that
$\mathcal{O}_1\cup\mathcal{O}_2\cup\mathcal{O}_4$ is a two-graph and (using
that $\tau$ is a proper two-graph) we have
$\tau=\mathcal{O}_1\cup\mathcal{O}_2\cup\mathcal{O}_4$. 

Recall that the Kneser graph $\Gamma$ is the graph with vertex set $\Omega$
and where two vertices $\alpha$ and $\beta$ are declared to be adjacent if 
they are disjoint, that is, $\alpha\cap\beta=\emptyset$. It is readily seen 
that the two-graph afforded by $\Gamma$ is $\tau$ and hence
$\Gamma\in \mathcal{S}$. Thus $\Aut(\Gamma)\leq \Aut(\tau)=G$ and hence
$\Aut(\Gamma)<G$ because we are dealing with Type~II. It is well-known (see
for example~\cite[Theorem and Table~II--VI]{LPS1}) that, for
$m\geq 9$, $\Aut(\Gamma)\cong S_m$ and $\Aut(\Gamma)$ is a
maximal subgroup of $\Sym(\Omega)$. Thus $\Aut(\Gamma)=\Aut(\mathcal{S})$,
contradicting the assumption of Type~II. 

When $m\leq 8$ and $k=2$, we have $n\in \{10,15,21,28\}$ and hence $n\leq 32$.

Suppose now that $k\geq 3$.  A careful analysis of the elements of $G$
yields $\fix_\Omega(g)\leq {m-2\choose k}+{m-2\choose k-2}$ for every 
$g\in G\setminus\{1\}$; the upper bound is achieved when $g$ is a
transposition. Observe also that $\orb_\Omega(g)\leq n/2$ for every $g\in G_2$.
Therefore from Equation~(\ref{eq:eq2}) and Lemma~\ref{basic3} we deduce
\begin{equation}
\label{eq:eq33}
2^{{m\choose k}}\leq 2^{\frac{1}{2}{m\choose k}+\frac{1}{2}\left({m-2\choose k}+{m-2\choose k-2}\right)}|\Aut(A_m)|
+2^{\frac{1}{2}{m\choose k}}|\Aut(A_m)|^2.
\end{equation}
Thus 
\[2^{{m\choose k}}\leq 2^{\frac{1}{2}{m\choose k}+\frac{1}{2}\left({m-2\choose k}+{m-2\choose k-2}\right)+1}|\Aut(A_m)|^2\]
and
\[2^{\frac{1}{2}{m\choose k}-\frac{1}{2}{m-2\choose k}-\frac{1}{2}{m-2\choose k-2}}\leq 2|\Aut(A_m)|^2,\]
and hence
\begin{eqnarray}\label{eq:eq4}
2^{{m-2\choose k-1}}\leq 2|\Aut(A_m)|^2,
\end{eqnarray}
where we use Pascal's recurrence for binomial coefficients twice to obtain
\[{m\choose k}={m-2\choose k}+2{m-2\choose k-1}+{m-2\choose k-2}.\]
Using $2^{{m-2\choose k-1}}\geq 2^{(m-2)(m-3)/2}$ and  $m!\leq m^{m-1}$,
we find that Equation~(\ref{eq:eq4}) holds true only for $m\leq 21$. 
For $m\leq 21$, a careful computation yields that Equation~(\ref{eq:eq33}) 
holds true only for $m\leq 8$. Finally, when $m\leq 8$, we may consider in
turn all the possibilities for $G$ and check that Equation~(\ref{eq:eq2}) is
never satisfied.

\subparagraph{Case (a)(ii)}
Now suppose that the socle of $G$ is $A_m^\ell$ in its 
product action on $\ell$ direct copies of $\{1,\ldots,m\}$ with $\ell\geq 2$
and $m\geq 5$.  

\smallskip

\noindent Suppose first that $\ell=2$. Hence $n=m^2$ and we may identify
$\Omega$ with the set of ordered pairs of elements  of $\{1,\ldots,m\}$.
Let $\tau$ be the two-graph induced by the switching class $\mathcal{S}$. 
Now, $G$ has rank $3$ in its action on $\Omega$ and  $G$ has four orbits on
the set of $3$-subsets of $\Omega$: namely
\begin{eqnarray*}
\mathcal{O}_1&=&\{\{(a,a),(a,b),(a,c)\}\mid a,b,c\textrm{ distinct elements of }\{1,\ldots,m\}\},\\
\mathcal{O}_2&=&\{\{(a,a),(a,b),(c,c)\}\mid a,b,c\textrm{ distinct elements of }\{1,\ldots,m\}\},\\
\mathcal{O}_3&=&\{\{(a,a),(a,b),(b,a)\}\mid a,b\textrm{ distinct elements of }\{1,\ldots,m\}\},\\
\mathcal{O}_4&=&\{\{(a,a),(b,b),(c,c)\}\mid a,b,c\textrm{ distinct elements of }\{1,\ldots,m\}\}.
\end{eqnarray*} 
(These orbits correspond to the non-equivalent geometric positions of three 
points in an $m\times m$ grid.) Replacing $\tau$ by its complement 
$\tau'=\{x\subseteq\Omega\mid |x|=3, x\notin \tau\}$, we may assume  that 
$\mathcal{O}_1\subseteq \tau$. Now consider the $4$-subset 
$x=\{(a,a),(a,b),(a,c),(b,a)\}\}$ of $\Omega$. Clearly, $x$ has one $3$-subset 
in common with $\mathcal{O}_1$, one with $\mathcal{O}_2$, two with 
$\mathcal{O}_3$ and zero with $\mathcal{O}_4$. As $\tau$ is a two-graph, 
$x$ has an even number of $3$-subsets in common with $\tau$ and hence 
$\mathcal{O}_2\subseteq \tau$. Now an easy inspection gives that 
$\mathcal{O}_1\cup\mathcal{O}_2$ is a two-graph and (using that $\tau$ is a 
proper two-graph) we have $\tau=\mathcal{O}_1\cup\mathcal{O}_2$. 

Recall that the grid graph $\Gamma'$ is the graph with vertex set $\Omega$ 
and where two distinct vertices $\alpha$ and $\beta$ are declared to be adjacent if 
$\alpha$ and $\beta$ are in the same row or in the same column (that is, the
first or the second coordinates of $\alpha$ and $\beta$ are equal). It is 
readily seen that the two-graph afforded by the complement $\Gamma$ of 
$\Gamma'$ is $\tau$ and hence $\Gamma\in \mathcal{S}$. Thus 
$\Aut(\Gamma)\leq \Aut(\tau)=G$ and hence $\Aut(\Gamma)<G$ because we are 
dealing with Type~II. It is well-known (see for 
example~\cite[Theorem and Table~I]{LPS1}) that 
$\Aut(\Gamma)\cong S_m\wr S_2$ and $\Aut(\Gamma)$ is a maximal 
subgroup of $\Sym(\Omega)$. Thus $\Aut(\Gamma)=\Aut(\mathcal{S})$, 
contradicting the assumption of Type~II.

Suppose then that $\ell\geq 3$. Now $\fix_\Omega(g)\leq m^{\ell-1}(m-2)$ for 
every $g\in G\setminus\{1\}$; the upper bound is achieved when $g$ is a 
transposition (see for example~\cite[Lemma~$6.13$]{GPS} or the proof 
of~\cite[Lemma~$4.5$]{Spiga1}). Observe also that $\orb_\Omega(g)\leq n/2$ for 
every $g\in G_2$. Therefore from Equation~(\ref{eq:eq2}) and 
Lemma~\ref{basic3} we deduce
\begin{equation}
\label{eq:eq44}
2^{m^\ell}\leq 2^{\frac{1}{2}{m^\ell}+
\frac{1}{2}\left(m^{\ell-1}(m-2)\right)}m!^\ell\ell!
+2^{\frac{1}{2}{m^\ell}}(m!^\ell\ell!)^2.
\end{equation}
Thus 
\[2^{m^\ell}\leq 2^{m^\ell-m^{\ell-1}+1}(m!^\ell\ell!)^2\]
and hence
\begin{equation}\label{eq:eq444}
2^{m^{\ell-1}-1}\leq (m!^\ell\ell!)^2.
\end{equation}
Using $\ell\geq 3$, $m!\leq m^{m-1}$ and $\ell!\leq \ell^{\ell-1}$, we find
that Equation~(\ref{eq:eq444}) holds true only when $\ell=3$. 
In this case, a careful computation yields that Equation~(\ref{eq:eq44}) is
never satisfied.

\subparagraph{Case (a)(iii):} the socle of $G$ is $A_m^\ell$ in its 
product action on $\ell$ direct copies of the $k$-subsets of $\{1,\ldots,m\}$
with $\ell\geq 2$, $m\geq 5$ and $2\leq k< m/2$.   

\smallskip

\noindent The argument is similar to the previous two cases. Here 
$n={m\choose k}^\ell$. From~\cite[Lemma~$6.13$]{GPS} (or from an easy 
computation), we deduce that
\[\fix_\Omega(g)\leq {m\choose k}^{\ell-1}\left({m-2\choose k}+{m-2\choose k-2}\right),\]
for every $g\in G\setminus\{1\}$. Moreover, as usual, 
$\fix_\Omega(g)\leq 2^{n/2}$ for every $g\in G_2$. Therefore from
Equation~(\ref{eq:eq2}) and Lemma~\ref{basic3} we get
\begin{equation}
\label{eq:eq55}
2^{{m\choose k}^\ell}\leq 2^{\frac{1}{2}{{m\choose k}^\ell}+
\frac{1}{2}\left({m\choose k}^{\ell-1}\left({m-2\choose k}+{m-2\choose k-2}\right)\right)}m!^\ell\ell!
+2^{\frac{1}{2}{{m\choose k}^\ell}}(m!^\ell\ell!)^2.
\end{equation}

Observe that 
\[\frac{1}{2}{{m\choose k}^\ell}
+\frac{1}{2}\left({m\choose k}^{\ell-1}\left({m-2\choose k}+{m-2\choose k-2}\right)\right)
={m\choose k}^\ell-{m-2\choose k-1}{m\choose k}^{\ell-1}.\]
Now, using $k\geq 2$, $\ell\geq 2$ and ${m\choose k}\geq m(m-1)/2$ and arguing
as in the previous two cases, we see that Equation~(\ref{eq:eq55}) is never
satisfied.

\subparagraph{Case (b):} the group $G$ is the Mathieu group $M_{n}$ with
$n\in \{11,12,23,24\}$ and hence $n\leq 32$. 



\subparagraph{Case (c):} We assume that none of the previous
cases occurs.

\smallskip

\noindent From~\cite[Theorem~$1.1$]{maroti}, we have $|G|\leq n^{1+\log_2(n)}$.
Moreover from~\cite[Corollary~$1$]{GM}, for every $g\in G\setminus\{1\}$, we
have $\fix_\Omega(g)\leq 4n/7$. In particular, for every
$g\in G\setminus\{1\}$, we get $\orb_\Omega(g)\leq n/2+2n/7= 11n/14$ by
Lemma~\ref{basic3}. Observe also that, for every $g\in G_2$,
$\orb_\Omega(g)\leq n/2$. Taking these bounds into account,
Equation~(\ref{eq:eq2}) yields
\[
\frac{2^{n}}{|G|}\leq \frac{1}2\left(\frac{2^n}{|G|}+2^{\frac{11n}{14}}\frac{|G|-1}{|G|}+|G_2|2^{\frac{n}{2}}\right)< \frac{1}{2}\left(\frac{2^{n}}{|G|}+2^{\frac{11n}{14}}+|G|2^{\frac{n}{2}}\right).
\]
Hence
\begin{equation}\label{eq:eq3}
2^{n}\leq 2^{11n/14}|G|+2^{n/2}|G|^2\leq 2^{11n/14}n^{1+\log_2(n)}+2^{n/2}n^{2+2\log_2(n)}.
\end{equation}
A computation yields $n\leq 384$. 

Observe that $G$ is not $3$-homogeneous because it is the 
automorphism group of a non-trivial two-graph. In particular, we may (and will) assume that $G$ is not $3$-homogeneous. Now that the degree $n$ is so 
small we can afford to compute the exact value of Equation~(\ref{eq:eq2}).
Indeed, a computer calculation using the database of small primitive groups
shows that, for $n\leq 384$
,  Equation~(\ref{eq:eq2}) holds
true only if $n\leq 64$.

For the remaining cases, we first construct the switching module $W$ and we
compute (using the Orbit-Stabiliser lemma) the exact value of
\[\#G\textrm{-orbits on }W=\frac{1}{|G|}\sum_{g\in G}|\cent W g|.\]
Using this formula and  Equation~(\ref{eq:eq0}), we check that
Equation~(\ref{hopefullylast})
is satisfied only for primitive groups of degree $n\leq 32$.
\end{pf}

\section{Groups of small degree}\label{sec5}

For groups of small degree, including the examples, we adopted a different
strategy. We only need to consider even degrees $n$: as explained earlier,
any possible example for odd $n$ will fall under Type~I, and will be on
Seress' list. A simplified version
of the algorithm takes a primitive group $G$ of degree $n$, and does the
following.
\begin{enumerate}
\item Compute the orbits of $G$ on $3$-element subsets of $\{1,\ldots,n\}$.
\item For each union of orbits, check whether $G$ is the full automorphism
group of the corresponding $3$-uniform hypergraph, and discard it if not.
Then check whether the hypergraph is a two-graph, and discard it if not.
\item On reaching this point, compute the graphs in the corresponding switching
class and their automorphism groups. Stop when either a graph with trivial
group is found, or every graph in the switching class has been considered.
In the latter case, record that an example has been found.
\end{enumerate}
The computations were done using \textsf{GAP}~\cite{gap}. We describe each
step in more detail.

Step (a) is straightforward; a single line of \textsf{GAP} code does this.
In Step  (b), we only need one of each complementary pair of two-graphs; this
is most easily achieved by omitting the last orbit on triples from the union.
Also, we do not have to consider the empty collection of orbits. We compute
the automorphism group using \texttt{nauty}~\cite{nauty}, interfaced to 
\textsf{GAP} via the DESIGN package~\cite{design}.

To check the two-graph property, it would suffice to check orbit
representatives for $G$ on $4$-subsets, to see whether each contains
an even number of $3$-subsets from the collection being considered. To avoid a
potentially large orbit computation, but use the fact that $G$ is transitive,
we simply checked all $4$-subsets containing the point $n$.

Up to degree $32$, the only two-graphs on an even number of points with
primitive automorphism groups are
\begin{itemize}
\item those with $2$-transitive groups, classified by Taylor~\cite{taylor}:
these are the Paley two-graphs with automorphism group 
$\mathrm{P}\Sigma\mathrm{L}(2,q)$, where $q$ is a prime power congruent to
$1$~(mod~$4$), the symplectic two-graph on $16$ points with group $2^4.S_6$,
and the orthogonal two-graph on $28$ points with group $\mathrm{PSp}(6,2)$; 
\item two (isomorphic) examples on $10$ points with group $A_5$, six on $28$
points with group $\mathrm{PGL}(2,7)$, and six on $28$ points with group
$\mathrm{PSL}(2,8)$.
\end{itemize}
Exploration in the range from $32$ to $40$ suggests that examples become
commoner.

The second type above surprised us a little; here is an explanation of the
two examples on $10$ points. The group is $A_5$ acting on pairs. The orbits
of $S_5$ on triples of pairs are isomorphism types of graphs with five vertices
and three edges, namely $K_3$, $K_{1,3}$, $K_2\cup P_3$, and $P_4$, where
$K_r$, $K_{r,s}$ and $P_s$ are complete graphs, complete bipartite graphs, and
paths respectively (the subscripts are the number of vertices). Since the
automorphism group of $P_4$ contains only even permutations, this orbit splits
into two under the action of $A_5$. Now Table~\ref{a5twographs} gives the
inclusions of graphs in these orbits in $4$-edge graphs.

\begin{table}[htb]
\begin{center}
\begin{tabular}{|c||c|c|c|c|c|}
\hline
&
\setlength{\unitlength}{0.5mm}
\begin{picture}(20,20)
\multiput(5,2)(10,0){2}{\circle*{1}}
\multiput(2,11)(16,0){2}{\circle*{1}}
\put(10,17){\circle*{1}}
\put(2,11){\line(1,0){16}}
\put(2,11){\line(4,3){8}}
\put(18,11){\line(-4,3){8}}
\end{picture}
&
\setlength{\unitlength}{0.5mm}
\begin{picture}(20,20)
\multiput(5,2)(10,0){2}{\circle*{1}}
\multiput(2,11)(16,0){2}{\circle*{1}}
\put(10,17){\circle*{1}}
\put(10,17){\line(-1,-3){5}}
\put(10,17){\line(1,-3){5}}
\put(10,17){\line(-4,-3){8}}
\end{picture}
&
\setlength{\unitlength}{0.5mm}
\begin{picture}(20,20)
\multiput(5,2)(10,0){2}{\circle*{1}}
\multiput(2,11)(16,0){2}{\circle*{1}}
\put(10,17){\circle*{1}}
\put(10,17){\line(-4,-3){8}}
\put(10,17){\line(4,-3){8}}
\put(5,2){\line(1,0){10}}
\end{picture}
&
\setlength{\unitlength}{0.5mm}
\begin{picture}(20,20)
\multiput(5,2)(10,0){2}{\circle*{1}}
\multiput(2,11)(16,0){2}{\circle*{1}}
\put(10,17){\circle*{1}}
\put(5,2){\line(-1,3){3}}
\put(10,17){\line(-4,-3){8}}
\put(10,17){\line(4,-3){8}}
\end{picture}
&
\setlength{\unitlength}{0.5mm}
\begin{picture}(20,20)
\multiput(5,2)(10,0){2}{\circle*{1}}
\multiput(2,11)(16,0){2}{\circle*{1}}
\put(10,17){\circle*{1}}
\put(5,2){\line(-1,3){3}}
\put(10,17){\line(-4,-3){8}}
\put(10,17){\line(4,-3){8}}
\end{picture}
\\
\hline\hline
\setlength{\unitlength}{0.5mm}
\begin{picture}(20,20)
\multiput(5,2)(10,0){2}{\circle*{1}}
\multiput(2,11)(16,0){2}{\circle*{1}}
\put(10,17){\circle*{1}}
\put(2,11){\line(1,0){16}}
\put(2,11){\line(4,3){8}}
\put(18,11){\line(-4,3){8}}
\put(5,2){\line(-1,3){3}}
\end{picture}
&1&1&0&1&1\\
\hline
\setlength{\unitlength}{0.5mm}
\begin{picture}(20,20)
\multiput(5,2)(10,0){2}{\circle*{1}}
\multiput(2,11)(16,0){2}{\circle*{1}}
\put(10,17){\circle*{1}}
\put(2,11){\line(1,0){16}}
\put(2,11){\line(4,3){8}}
\put(18,11){\line(-4,3){8}}
\put(5,2){\line(1,0){10}}
\end{picture}
&1&0&3&0&0\\
\hline
\setlength{\unitlength}{0.5mm}
\begin{picture}(20,20)
\multiput(5,2)(10,0){2}{\circle*{1}}
\multiput(2,11)(16,0){2}{\circle*{1}}
\put(10,17){\circle*{1}}
\put(5,2){\line(1,0){10}}
\put(2,11){\line(1,0){16}}
\put(5,2){\line(-1,3){3}}
\put(15,2){\line(1,3){3}}
\end{picture}
&0&0&0&2&2\\
\hline
\setlength{\unitlength}{0.5mm}
\begin{picture}(20,20)
\multiput(5,2)(10,0){2}{\circle*{1}}
\multiput(2,11)(16,0){2}{\circle*{1}}
\put(10,17){\circle*{1}}
\put(5,2){\line(1,3){5}}
\put(2,11){\line(4,3){8}}
\put(18,11){\line(-4,3){8}}
\put(15,2){\line(1,3){3}}
\end{picture}
&0&1&1&1&1\\
\hline
\setlength{\unitlength}{0.5mm}
\begin{picture}(20,20)
\multiput(5,2)(10,0){2}{\circle*{1}}
\multiput(2,11)(16,0){2}{\circle*{1}}
\put(10,17){\circle*{1}}
\put(10,17){\line(-4,-3){8}}
\put(10,17){\line(4,-3){8}}
\put(10,17){\line(-1,-3){5}}
\put(10,17){\line(1,-3){5}}
\end{picture}
&0&4&0&0&0\\
\hline
\setlength{\unitlength}{0.5mm}
\begin{picture}(20,20)
\multiput(5,2)(10,0){2}{\circle*{1}}
\multiput(2,11)(16,0){2}{\circle*{1}}
\put(10,17){\circle*{1}}
\put(5,2){\line(-1,3){3}}
\put(15,2){\line(1,3){3}}
\put(10,17){\line(-4,-3){8}}
\put(10,17){\line(4,-3){8}}
\end{picture}
&0&0&2&0 or 2&2 or 0\\
\hline
\end{tabular}
\end{center}
\caption{\label{a5twographs} Graphs on $5$ vertices}
\end{table}

The table shows (up to complementation) one two-graph admitting $S_5$,
consisting of the $P_4$ graphs. This is the Paley two-graph with
automorphism group $\mathrm{P}\Sigma\mathrm{L}(2,9)$, aka $S_6$. But we have
two further two-graphs admitting $A_5$, each consisting of one orbit on $P_4$s
together with the $K_{1,3}$s. An element of $S_5\setminus A_5$ induces an
isomorphism between these two-graphs; their symmetric difference is the Paley
two-graph.

For step (c), we first construct the graph in the switching class which has
vertex $n$ isolated: join $x$ and $y$ if and only if $\{x,y,n\}$ is a triple
in the collection. We can then use the fact that switching with respect to
a set and its complement are identical, so we need only consider switching
sets not containing $n$. These sets are generated iteratively, the starting
graph is switched, and the automorphism group found by calling \texttt{nauty}
from the \textsf{GAP} package GRAPE~\cite{grape}.

\end{document}